\documentclass[letterpaper]{article}

\usepackage[margin=2.5cm]{geometry}
\usepackage{amsmath,amsthm,amsfonts,amssymb, array, tabularx, url}
\usepackage[sortcites=true]{biblatex}
\newcommand{\overbar}[1]{\mkern 1.5mu\overline{\mkern-1.5mu#1\mkern-1.5mu}\mkern 1.5mu}

\newcommand{\E}{\mathbb E}
\newcommand{\bbF}{\mathbb F}

\newcommand{\ex}{\mathrm{ex}}
\newcommand{\CD}{\mathrm{CD}}
\newcommand{\D}{\mathrm{D}}
\newcommand{\Aut}{\mathrm{Aut}}
\newcommand{\PG}{\mathrm{PG}}

\newtheorem{thm}{Theorem}[section]
\newtheorem{lemma}[thm]{Lemma}
\newtheorem{prop}[thm]{Proposition}
\newtheorem{prob}[thm]{Problem}

\theoremstyle{definition}
\newtheorem{cor}[thm]{Corollary}
\newtheorem{clm}[thm]{Claim}
\newtheorem{defn}[thm]{Definition}

\newtheorem*{rmk}{Remark}
\numberwithin{equation}{section}
\numberwithin{figure}{section}

\title{Extremal Numbers of Hypergraph Suspensions of Even Cycles}

\author{
	Sayan Mukherjee\footnote{Blueqat Research, Tokyo 150-6139, Japan.} \footnote{Department of Physics, The University of Tokyo, Tokyo 113-0033, Japan. Email: \texttt{sayan@g.ecc.u-tokyo.ac.jp}. \\ \hspace*{0.5cm}(Article accepted at European Journal of Combinatorics)}
}

\begin{document}
	\maketitle
    \vspace*{-0.8cm}
	\begin{abstract}
		For fixed $k\ge 2$, determining the order of magnitude of the number of edges in an $n$-vertex bipartite graph not containing $C_{2k}$, the cycle of length $2k$, is a long-standing open problem.
		We consider an extension of this problem to triple systems.
		In particular, we prove that the maximum number of triples in an $n$-vertex triple system which does not contain a $C_6$ in the link of any vertex, has order of magnitude $n^{7/3}$.
		Additionally, we construct new families of dense $C_6$-free bipartite graphs with $n$ vertices and $n^{4/3}$ edges in order of magnitude.
	\end{abstract}
\vspace*{-0.5cm}
\begin{center}{\small\textbf{Keywords:} Hypergraph Tur\'an problem, algebraic constructions, hypergraph suspension, even cycles.}\end{center}
\vspace*{-0.5cm}
\section{Introduction}

An $r$-uniform hypergraph, or simply, an $r$-graph $H$ on vertex set $V(H)$ is a subset of $\binom{V(H)}{r}$.
In this work, we denote by $|H|$ the number of edges, or simply, the size of $H$.
For a fixed $r$-graph $H$, we say that an $r$-graph $G$ is \emph{$H$-free} if $G$ does not contain a copy of $H$ as a subgraph.
The hypergraph Tur\'an problem asks the following question: what is the largest number of edges in an $H$-free $r$-graph on $n$ vertices?
This number is known as the Tur\'an number or the extremal number of $H$, and is denoted by $\ex_r(n,H)$.
The case $r=2$ was first introduced by Tur\'an \cite{turanExtremal1941} in 1941, and several lower and upper bounds on $\ex_r(n,H)$ have been obtained since then for different values of $r$ and $H$.

Towards analyzing the asymptotic behavior of $\ex_2(n,G)$ for graphs $G$, the seminal result of Erd\H os and Stone \cite{erdosStoneLinear1963} states that when the chromatic number $\chi(G)$ is at least $3$,
\[
\ex_2(n,G) = \left(1-\frac1{\chi(G)-1}\right)\binom n2 + o(n^2).
\]
This result essentially determines $\ex_2(n,G)$ for graphs $G$ which are not bipartite.
The analysis of $\ex_2(n,G)$ for bipartite graphs $G$ turns out to be extremely difficult, and the reader is referred to \cite{furediSimonovitsDegenerateSurvey2013} for a comprehensive survey of the bipartite case.

One especially well-studied class of bipartite graphs $G$ are the even cycles $C_{2k}$ for $k\ge 2$.
For these graphs, the best known general upper bound is due to Bondy and Simonovits \cite{bondySimonovitsEvenCycle1974}, who proved that $\ex_2(n,C_{2k}) = O(n^{1+\frac1k})$.
Improvements in the coefficient of $n^{1+\frac 1k}$ has been obtained in \cite{bollobasExtremalGraphTheory1978, verstraeteArithmeticProgression2000, PikhurkoEvenCycle2012, bukhJiangEvenCycle2017}.

A major open problem for even cycles is to construct $C_{2k}$-free graphs on $n$ vertices and $\Omega(n^{1+\frac1k})$ edges.
There have been several bipartite constructions based on finite fields \cite{reimanProblem1957, brownThomsen1966, bensonminimal1966, margulisExplicit1982, wengerExtremal1991, lazebnikUstimenkoWoldarPolarity1999, terlepWillifordKacMoody2012} that have provided lower bounds on specific values of $k$.
For general $k$ (except $k\in \{2,3,5,7\}$) the best known lower bounds are given by the bipartite graphs $\CD(k,q)$ for integers $k\ge 2$ and prime powers $q$ \cite{lazebnikGirth1995, lazebnikUstimenko1995}.
These graphs arise from Lie algebraic incidence structures that approximate the behavior of generalized polygons, and are analyzed in detail in \cite{woldarGeneralizedPoly2009}.
See also \cite{lazebnikViglioneConnectivity2004} for an analysis of the connectivity of $\CD(k,q)$ for even $q$. 
For a recent survey on the even cycle problem, the reader is referred to \cite{verstraeteCycleSurvey2016}.

In this paper, we are mainly concerned with three classes of lower bound constructions: the bipartite graphs $\D(k,q)$ from \cite{lazebnikGirth1995}, the {arc construction} introduced in \cite{mellingerLDPC2004} and later generalized in \cite{mellingerMubayiArc2005}, and Wenger's construction~\cite{wengerExtremal1991}.
Our results can be divided into two sections: results about 3-graphs and results about graphs.

\subsection{3-Graphs}
For a graph $G$, the suspension $\widehat G$ is the graph obtained from $G$ by adding a new vertex adjacent to all vertices of $G$.
For fixed $n$, the \emph{generalized Tur\'an number} $\ex(n,T,H)$ (studied rigorously in \cite{alonTCopies2016}) is defined as the maximum number of (non-induced) copies of $T$ in an $H$-free graph on $n$ vertices.
In~\cite{suspensionFree2020}, the author, together with Mubayi, studied $\ex(n, K_3, \widehat G)$ for different bipartite graphs $G$.
Analogously, we introduce the concept of a \emph{hypergraph suspension}.

Let $H$ be a 3-graph and $x\in V(H)$ be any vertex of $H$.
The link of $x$ in $H$, denoted by $L_{x,H}$, is the graph with vertex set $V(H)\setminus \{x\}$ and edges $\{uv: \{x,u,v\}\in H\}$.
For a graph $G$, the hypergraph suspension $\widetilde{G}$ is a 3-graph defined as follows: add a new vertex $x$ to $V(G)$, and let $\widetilde G=\{e\cup\{x\}:e\in E(H)\}$.
By definition, $L_{x,\widetilde G} = G$.

Note that the numbers $\ex_3(n,\widetilde G)$ and $\ex(n,K_3,\widehat G)$ are closely related.
In fact, given a $\widehat G$-free graph, we can replace all triangles in it with hyperedges to obtain a $\widetilde G$-free 3-graph, implying
\begin{equation}
\label{eq:hat-tilde-reln}
\ex(n,K_3,\widehat G)\le \ex_3(n,\widetilde{G}).
\end{equation}

In this paper, we study $\ex_3(n,\widetilde C_{2k})$ for $k\ge 2$.
When $k=2$, observe that $\widetilde C_{2k}$ is the complete 3-partite 3-graph $K^{(3)}_{1,2,2}$, and thus it is shown in \cite{mubayiTripleC4-2002} that $\ex_3(n,\widetilde C_4)=\Theta(n^{5/2})$.
Hence we focus our attention on $\ex_3(n,\widetilde C_{2k})$ for $k\ge 3$.

Observe that a $3$-graph $H$ does not contain $\widetilde C_{2k}$ iff $L_{x,H}$ does not contain $C_{2k}$ for every vertex $x\in V(H)$, leading us to the upper bound
\begin{equation}
\label{eq:ex3-n-C2k-upper}
\ex_3(n,\widetilde C_{2k})= O(n\cdot n^{1+\frac 1k}) = O(n^{2+\frac 1k})
\end{equation}

On the other hand, a probabilistic deletion argument lets us deduce the following result:

\begin{prop}
	\label{prop:probabilistic-lower-C2k}
	For $k\ge 2$,
	\begin{equation}
	\label{eq:ex3-n-C2k-lower-prob}
	ex_3(n,\widetilde C_{2k})= \Omega\left(n^{2+\frac 1{2k-1}}\right).
	\end{equation}
\end{prop} 

Our main result is to show a construction of $\widetilde C_{2k}-$free 3-graphs, which asymptotically improves the bound above for $k=3$ and $k=4$.

\begin{thm}
	\label{thm:c2ktilde}
	For every integer $q$ that is a power of $3$, there exists a 3-partite 3-graph $\D_3(k,q)$ with the following properties:
	\begin{enumerate}
	    \item $\D_3(k,q)$ has $3q^k$ vertices and $q^{2k+1}$ edges,
	    \item The link graphs of every vertex of $\D_3(k,q)$ are isomorphic for $k\le 6$, and
	    \item $\D_3(3,q)$ and $\D_3(5,q)$ are $\widetilde C_6$-free and $\widetilde C_8$-free, respectively.
	\end{enumerate}
\end{thm}

In particular, Theorem~\ref{thm:c2ktilde} implies that\footnote{
Since $\ex_3(n,H)$ is monotonically increasing in $n$, a lower bound of $\ex_3(3^r,H)\ge C\cdot 3^{r\alpha}$ implies $\ex_3(n,H)\ge \ex_3(3^{\lfloor\log_3 n\rfloor}, H) \ge C\cdot 3^{\lfloor \log_3 n\rfloor\alpha} \ge 3^{-\alpha} C\cdot n^\alpha$. Thus an asymptotic lower bound on powers of $3$ easily generalizes to all natural numbers $n$.
}
\begin{equation}
\label{eq:ex3-n-C2k-lower}
    \ex_3(n,\widetilde C_{6}) = \Omega(n^{7/3})\mbox{ and }\ex_3(n,\widetilde C_{8}) = \Omega(n^{11/5}).
\end{equation}

As a corollary of (\ref{eq:ex3-n-C2k-upper}) and (\ref{eq:ex3-n-C2k-lower}), we determine the asymptotic growth rate of $\ex_3(n,\widetilde C_6)$.
\begin{cor}
\label{cor:ex3-n-C6-exact}
For large $n$, the Tur\'an number of $\widetilde C_6$ grows as,
\begin{equation}
    \ex_3(n,\widetilde C_{6}) = \Theta(n^{7/3}).
\end{equation}
\end{cor}
Corollary~\ref{cor:ex3-n-C6-exact} further implies that the bound in (\ref{eq:hat-tilde-reln}) is not always sharp, since we demonstrated in~\cite{suspensionFree2020} that $\ex(n,K_3,\widehat C_{6}) = o(n^{7/3})$.

\begin{rmk}
    Our proof of Theorem~\ref{thm:c2ktilde} heavily relies on the bipartite graphs $\D(k,q)$ introduced by Lazebnik, Ustimenko and Woldar \cite{lazebnikGirth1995}, and $\D_3(k,q)$ can be viewed as an extension of $\D(k,q)$ to 3-graphs. $\D_3(k,q)$ has the property that for every $k\ge 2$ and $q$ a power of $3$, the link graph of any of its vertex is isomorphic to either $\D(k,q)$ or another graph which we call $\D'(k,q)$ (Proposition~\ref{prop:D3kq-twolinks}). We also propose a problem of figuring out the girth and connectivity of $\D'(k,q)$ (Problem~\ref{prob:girthdpkq}). A good lower bound on the girth of  $\D'(k,q)$ would directly translate to asymptotic improvements in the lower bound on $\ex_3(n,\widetilde C_{2k})$ given by (\ref{eq:ex3-n-C2k-lower-prob}).
\end{rmk}

\subsection{Graphs}
We also compare two well-known constructions of $C_{2k}$-free graphs: the arc construction~\cite{mellingerLDPC2004,mellingerMubayiArc2005} and Wenger's construction~\cite{wengerExtremal1991}.
Before describing these constructions, we introduce some basic notation from (finite) projective geometry.

Let $t\ge 2$, and let $q$ be a prime power.
Let $\bbF_q$ denote the finite field of $q$ elements.
The vector space $V=\bbF_q^{t+1}$ naturally defines a $(t+1)$-dimensional vector space over $\bbF_q$.
Let $\sim$ be the equivalence relation on $V$ defined by $x\sim y$ iff there exists $\lambda\in \bbF_q$, $\lambda\neq 0$ such that $x = \lambda y$.
Then, the projective space $\PG(t,q)$ is the set of equivalence classes of $V \setminus \{0\}$ under $\sim$.
Note that since $\dim_{\bbF_q} V = t+1$, we have $\dim_{\bbF_q}(\PG(t,q)) = t$.
Further, any point $x\in \PG(t,q)$ is usually given by its \emph{homogeneous coordinates}, $[x_0:x_1:\cdots:x_t]$, so that for any $\lambda\in \bbF_q$ with $\lambda\neq 0$, $[x_0:x_1:\cdots:x_t]=[\lambda x_0:\lambda x_1: \cdots: \lambda x_t]$.
An arc in a projective $t$-space $\PG(t,q)$ is a collection of points such that no $(t-1)$ of them lie in a hyperplane.

Now we present the defition of the arc construction~\cite{mellingerLDPC2004,mellingerMubayiArc2005}.

\textbf{The bipartite graphs $G_{\texttt{arc}}(k,q,\alpha)$.}
Let $\Sigma = \PG(t,q)$, and $\Sigma_0\subset \Sigma$ be the hyperplane consisting of points with first homogeneous coordinate $0$.
Note that $\Sigma_0\cong \PG(t-1,q)$.
Let $\alpha$ be any arc in $\Sigma_0$. Then, the bipartite graph $G_{\texttt{arc}}(k,q,\alpha)$ with parts $P$ and $L$ is defined as follows.
Let $P=\Sigma\setminus\Sigma_0$, and $L$ be the set of all projective lines $\ell$ of $\Sigma$ such that $\ell\cap\Sigma_0\in \alpha$.
Vertices $p\in P$ and $\ell\in L$ are adjacent if and only if $p\in \ell$.

It was shown by Mellinger and Mubayi \cite{mellingerMubayiArc2005} that $G_{\texttt{arc}}(k,q,\alpha_0)$ is $C_{2k}$-free for $k=2,3,5$ but contains $C_{2k}$ for $k=4$, where $\alpha_0$ is the normal rational curve in $\Sigma_0$ given by
\[
\alpha_0 = \{[0:1:x:x^2:\cdots:x^{t-1}]: x\in\bbF_q\}\cup\{[0:0:\cdots: 0: 1]\}.
\]

\textbf{The Wenger graphs $H(k,q)$.}
Let $H(k,q)$ be the bipartite graph with parts $A=B=\bbF_q^k$ such that $(a_1,\ldots, a_k)$ is adjacent to $(b_1,\ldots, b_k)$ iff \[
a_i+b_i = a_1 b_1^{i-1} \mbox{ for all } 2\le i \le k.
\]
It was shown by Wenger in~\cite{wengerExtremal1991} that $H(k,q)$ is $C_{2k}$-free for $k=2,3,5$.

\medskip
Cioab{\u{a}}, Lazebnik and Li \cite{cioabuaLazebnikLiSpectrumWenger2014} proved, among other properties of $H(k,q)$, that these two constructions are, in fact, isomorphic:

\begin{prop}
	\label{prop:Arc_Wenger_isomorphism}
	Let $\alpha_0$ be the normal rational curve in $\PG(k,q)$, and $\alpha_0^-=\alpha_0 \setminus\{[0:\cdots:0:1]\}$. Then,
	\[
	G_{\texttt{arc}}(k,q,\alpha_0^-)\cong H(k,q).
	\]
\end{prop}

We give an alternative proof of Proposition~\ref{prop:Arc_Wenger_isomorphism} using the Pl\"ucker embedding~\cite{griffithsHarrisAG1978}, a tool from algebraic geometry that lets us parametrize the set of projective lines $L$.


Let $(a,b)$ denote the greatest common divisor of integers $a$ and $b$.
For $1\le s \le r$ with $(s,r)=1$, it can be shown that (see, for example, Section 10.3.1 of \cite{payneThasFiniteGeneralizedQuadrangles2009} and Claim~\ref{clm:M(x,y,z)Invertible}) $\alpha = \{ [1:x:x^{2^s}]: x\in \bbF_{2^r} \}$ is an arc in the projective space $\PG(2,2^r)$.
Using the proof method of Proposition~\ref{prop:Arc_Wenger_isomorphism} on this arc $\alpha$, we are able to construct a family of $C_6$-free graphs with $\Omega(n^{4/3})$ edges, given as follows.
\begin{thm}
\label{thm:newC6construction}
	Let $q=2^r$ and $1\le s\le r$ be such that $(s,r)=1$. Let $G(2^r,s)$ denote the bipartite graph with parts $A=B=\bbF_q^3$ such that $(a_1, a_2, a_3)\in A$ is adjacent to $(b_1,b_2,b_3)\in B$ iff
	\[
	b_2+a_2 = a_1b_1 \mbox{ and } b_3+a_3 = a_1b_1^{2^s}.
	\]
	Then, $G(2^r, s)$ is $C_6$-free.
\end{thm}
Note that the graphs $G(2^r,s)$ extend Wenger's $C_6$-free construction in even characteristic, as $G(2^r,1) \cong H(3,2^r)$.
We are unsure about whether the graphs $G(2^r,s)$ are isomorphic to any other $C_6$-free graph families already known, in particular, whether $G(2^r,s)\cong G(2^r,1)$.
However, it might be possible to give explicit definitions for new bipartite constructions without even cycles of certain lengths using arcs in different projective spaces.

\medskip
This paper is organized as follows.
In Section 2, we prove Proposition~\ref{prop:probabilistic-lower-C2k}, recapitulate on the graphs $\D(k,q)$, extend them to the 3-graphs $\D_3(k,q)$, and investigate its link graphs, finally proving Theorem~\ref{thm:c2ktilde}.
Section 3 is devoted to proving Proposition~\ref{prop:Arc_Wenger_isomorphism} and Theorem~\ref{thm:newC6construction}.

\section{Lower bounds on $\ex_3(n,\widetilde C_{2k})$}
Our goal in this section is to extend the graphs $\D(k,q)$ to a family of 3-graphs, and build up the tools required to prove Theorem~\ref{thm:c2ktilde}. We start with a proof of Proposition~\ref{prop:probabilistic-lower-C2k}. Recall that we wish to show $\ex_3(n,\widetilde C_{2k})\ge \Omega(n^{2+\frac 1{2k-1}})$.

\begin{proof}[Proof of Proposition~\ref{prop:probabilistic-lower-C2k}]

Let $H\sim G_3(n,p)$ be the Erd\H os-R\'enyi $3$-graph, where each edge of the complete $3$-graph on $n$ vertices is selected independently with probability $p=c_kn^{-\frac{2k-2}{2k-1}}$ for a constant $c_k$ which we choose later. Then, $\E(|H|)=p\binom n3$. For every $\widetilde C_{2k}$ in $H$, we remove one edge from it. Let $H'\subset H$ be the new 3-graph obtained via the deletion of edges. Note that the probability that any $2k+1$ vertices forms a $\widetilde C_{2k}$ is $(2k+1)\cdot (2k)!/4k\cdot p^{2k}$, and therefore, the expected number of them is at most $(2k+1)!/4k \cdot n^{2k+1}p^{2k}$. Now, $\E(|H'|)=p\binom n3 - {(2k+1)!n^{2k+1}p^{2k}}/{4k}$. As
\[
n^{2k+1} p^{2k-1} = n^{2k+1} \cdot c_k^{2k-1}n^{-(2k-2)} = c_k^{2k-1}{n^3},
\]
we have
\[
\E(|H'|) = p\left(\binom n3-\frac{(2k+1)!n^{2k+1}p^{2k-1}}{4k}\right) \ge pn^3\left(\frac 1{10} - \frac{(2k+1)!c_k^{2k-1}}{4k}\right).
\]
Taking $c_k = \left(\frac{1}{100}\cdot \frac{4k}{(2k+1)!}\right)^{1/(2k-1)}$, we note that $\E(|H'|)\ge pn^3/100\ge \frac{c_k}{100}\cdot n^{3-\frac{2k-2}{2k-1}}$. Thus, there exists a $3$-graph $H'$ with $\Omega(n^{3-\frac{2k-2}{2k-1}})$ edges with no copy of $\widetilde C_{2k}$. This completes our proof.
\end{proof}

Since probabilistic lower bounds for 3-graphs tend to be weak, we try to strengthen this result via a look at the graphs $\D(k,q)$. Here we present a summary of the properties of $\D(k,q)$; for more details, the reader is referred to  \cite{lazebnikGirth1995, lazebnikUstimenko1995, woldarGeneralizedPoly2009}.

\begin{defn}[The bipartite graphs $\D(q)$]
	For a prime power $q$, let $A$ and $B$ be two \textit{disjoint} copies of the countably infinite dimensional vector space $V$ over $\bbF_q$. Use the following coordinate representations for elements $a\in A$ and $b\in B$:
	\begin{equation}
	\label{eq:coordinatesDq}
	\begin{array}{lllcclcll}
	a&=&(&\overbrace{a_1,a_{11},a_{12},a_{21}}^{\text{group }1}, &\overbrace{a_{22},a'_{22},a_{23},a_{32}}^{\text{group }2},&\ldots&\overbrace{a_{ii},a'_{ii},a_{i,i+1},a_{i+1,i}}^{\text{group }i},&\ldots&),\\
	b&=&(&b_1,b_{11},b_{12}, b_{21}, &{b_{22},b'_{22},b_{23},b_{32}},&\ldots &b_{ii},b'_{ii},b_{i,i+1},b_{i+1,i},&\ldots&).
	\end{array}
	\end{equation}
	Let $A\sqcup B$ be the vertex set of $\D(q)$, and join $a\in A$ to $b\in B$ if the following coordinate relations hold ($i\ge 2$):
	\begin{equation}
	\label{eq:defDq-ab}
    \mathcal R(a,b) :=\left\{
		\begin{array}{ccc}
		a_{11}+b_{11}+a_1b_1&=&0\\
		a_{12}+b_{12}+a_1b_{11}&=&0\\
		a_{21}+b_{21}+a_{11}b_1&=&0\\
		&\vdots&\\
		a_{ii}+b_{ii}+a_{i-1,i}b_1&=&0\\
		a'_{ii}+b'_{ii}+a_1b_{i,i-1}&=&0\\
		a_{i,i+1}+b_{i,i+1}+a_1b_{ii}&=&0\\
		a_{i+1,i}+b_{i+1,i}+a'_{ii}b_1&=&0\\
        &\vdots&
		\end{array}\right.
	\end{equation}
    Note that the first $k-1$ equations of (\ref{eq:defDq-ab}) has the first $k$ coordinates of (\ref{eq:coordinatesDq}).
    Let $A_k$ and $B_k$ denote the truncation of $A$ and $B$ from (\ref{eq:coordinatesDq}) to the first $k$ coordinates, and $\mathcal R_k$ the truncation of $\mathcal R$ from (\ref{eq:defDq-ab}) to the first $k$ relations.
    Then, $\D(k,q)$ is defined as the bipartite graph with bipartition $A_k\sqcup B_k$ where vertices $a\in A_k$ and $b\in B_k$ are adjacent if they satisfy $\mathcal R_{k-1}(a,b)$.
    
    Observe that for a fixed vertex $a\in A_k$ in $\D(k,q)$, the subspace $\{b\in B: \mathcal R_{k-1}(a,b)\text{ holds}\}$ has dimension $k-(k-1)=1$, implying that every $a\in A_k$ has $q$ neighbors in $B_k$.
	By symmetry, this is true for every vertex of $B_k$ as well, implying that $\D(k,q)$ is a $q$-regular graph on $2q^k$ vertices.
\end{defn}

The key properties of the graphs $\D(k,q)$ are summarized in the following proposition.

\begin{prop}
	\label{prop:dkqproperties}
	For any prime power $q$ and $k\ge 2$, the girth of $\D(k,q)$ is at least $k+4$ if $k$ is even, and $k+5$ if $k$ is odd.
\end{prop}

Further, it is known that for $k\ge 6$ the graphs $\D(k,q)$ start to get disconnected into pairwise isomorphic components at regular intervals.
These connected components are called $\CD(k,q)$.
The graphs $\CD(2k-3,q)$ give the currently best known asymptotic lower bounds on $\ex(n,C_{2k})$ for $k\ge 3$.
We omit the proof of Proposition~\ref{prop:dkqproperties} here.

In the following subsection, we extend $\D(k,q)$ to the 3-graph case.
\subsection{The 3-graphs $\D_3(k,q)$}


\begin{defn}[The 3-partite 3-graphs $\D_3(q)$]
	For a prime power $q$, let $A$, $B$, and $C$ be three \textit{disjoint} copies of the countably infinite dimensional vector space $V$ over $\bbF_q$. We use the following coordinate representations for $a\in A$, $b\in B$, $c\in C$:
	\begin{equation}
	\label{eq:coordinatesD3q}
	\begin{array}{lllcclcll}
	a&=&(&\overbrace{a_1,a_{11},a_{12},a_{21}}^{\text{group }1}, &\overbrace{a_{22},a'_{22},a_{23},a_{32}}^{\text{group }2},&\ldots&\overbrace{a_{ii},a'_{ii},a_{i,i+1},a_{i+1,i}}^{\text{group }i},&\ldots&),\\
	b&=&(&b_1,b_{11},b_{12}, b_{21}, &{b_{22},b'_{22},b_{23},b_{32}},&\ldots &b_{ii},b'_{ii},b_{i,i+1},b_{i+1,i},&\ldots&),\\
    c&=&(&c_1,c_{11},c_{12}, c_{21}, &{c_{22},c'_{22},c_{23},c_{32}},&\ldots &c_{ii},c'_{ii},c_{i,i+1},c_{i+1,i},&\ldots&).
	\end{array}
	\end{equation}
	Let $A\sqcup B\sqcup C$ be the vertex set of $\D_3(q)$, and say that $\{a,b,c\}$ is a hyperedge if the following coordinate relations hold ($i\ge 2$):
	\begin{equation}
	\label{eq:d3qDefn}
    \mathcal R^{(3)}(a,b,c):=\left\{
	\begin{array}{ccc}
	a_{11}+b_{11}+c_{11}+a_1b_1+b_1c_1+c_1a_1&=&0\\
	a_{12}+b_{12}+c_{12}+a_1b_{11}+b_1c_{11}+c_1a_{11}&=&0\\
	a_{21}+b_{21}+c_{21}+a_{11}b_1+b_{11}c_1+c_{11}a_1&=&0\\
	&\vdots&\\
	a_{ii}+b_{ii}+c_{ii}+a_{i-1,i}b_1+b_{i-1,i}c_1+c_{i-1,i}a_1&=&0\\
	a'_{ii}+b'_{ii}+c'_{ii}+a_1b_{i,i-1}+b_1c_{i,i-1}+c_1a_{i,i-1}&=&0\\
	a_{i,i+1}+b_{i,i+1}+c_{i,i+1}+a_1b_{ii}+b_1c_{ii}+c_1a_{ii}&=&0\\
	a_{i+1,i}+b_{i+1,i}+c_{i+1,i}+a'_{ii}b_1+b'_{ii}c_1+c'_{ii}a_1&=&0\\
    &\vdots&\\
	\end{array}\right.
	\end{equation}
    Let $A_k$, $B_k$, $C_k$ denote the truncations of $A$, $B$ and $C$ from (\ref{eq:coordinatesD3q}) to the first $k$ coordinates, and $\mathcal R^{(3)}_k$ the truncation of $\mathcal R^{(3)}$ from (\ref{eq:d3qDefn}) to the first $k$ relations.
    Define $\D_3(k,q)$ to be the 3-graph with vertex set $A_k\sqcup B_k\sqcup C_k$, such that $\{a,b,c\}$ is a hyperedge of $\D_3(k,q)$ if $\mathcal R^{(3)}_{k-1}(a,b,c)$ holds.
\end{defn}

For any vector $\vec v\in V$, let $\vec v_A\in A$, $\vec v_B\in B$ and $\vec v_C\in C$ denote the corresponding vertices of $\D_3(q)$.
We have designed the $3$-graphs $\D_3(q)$ in such a way that the equations governing the link graph of $\vec 0_A$, $\vec 0_B$, $\vec 0_C$ are the same as the equations defining $\D(q)$.

In fact, note that $\D_3(q)$ has the natural cyclic automorphism $a_\ast \mapsto b_\ast$, $b_\ast \mapsto c_\ast$, and $c_\ast \mapsto a_\ast$, under which all the defining equations of $\D_3(q)$ remain invariant.
Hence, for any $\vec v\in V$, the link graphs of $\vec v_A$, $\vec v_B$ and $\vec v_C$ are all isomorphic.
One would hope that the link graphs of vertices of $\D_3(k,q)$ corresponding to other vectors $\vec v\neq \vec 0$ would also have similar high girth properties as $\D(k,q)$.
This inspires us to analyze the links of every vertex in $\D_3(k,q)$.
To that end, we analyze $\Aut(\D_3(q))$.

\begin{prop}
	\label{prop:a-to-a1-000}
	Suppose $\bbF_q$ has characteristic $3$, and consider $\D_3(q)$ with parts $A$, $B$, $C$. Let $a\in A$ be fixed, and suppose $s\ge 1$. Then there is an automorphism $\varphi\in\Aut(\D_3(q))$ such that 
    \[\varphi(a)=(a_1,\overbrace{0,\ldots,0}^{s\text{ zeros}},\ast,\ast,\ldots)_A.\]
\end{prop}

The proof of Proposition~\ref{prop:a-to-a1-000} is technical.
Before looking at the proof, we note an important consequence:
to analyze the girths of every vertex of $\D_3(k,q)$, it is sufficient to analyze the girths of the link graphs of the vertices $(a_1,{0,\ldots,0})_{A_k}$ for $a_1\in\bbF_q$.
In fact, it is seen that the truncated 3-graphs $\D_3(k,q)$ have exactly two kinds of links.

\begin{prop}
    \label{prop:D3kq-twolinks}
    If $q$ is a power of $3$, then the 3-graph $\D_3(k,q)$ admits exactly two classes of link graphs, one of which is $\D(k,q)$.
\end{prop}

Now, we present the proofs of Propositions~\ref{prop:a-to-a1-000} and \ref{prop:D3kq-twolinks}.

\subsubsection{Proof of Proposition~\ref{prop:a-to-a1-000}}
Recall that $q$ is a power of $3$, and we wish to construct an automorphism $\varphi$ of $\D_3(q)$ sending any vertex $a\in A$ to \[(a_1,\overbrace{0,\ldots,0}^{s\text{ zeros}},\ast,\ast,\ldots)_A.\]

	We construct $\varphi$ via a product of automorphisms of $\D_3(q)$. First, we may rewrite the relations $\mathcal R^{(3)}$ from (\ref{eq:d3qDefn}) into the following form:
	\begin{equation}
	\label{eq:d3qDefnNew}
	\mathcal R^{(3)}(a,b,c) = \left\{\begin{aligned}
	a_{ii}+b_{ii}+c_{ii}+a_{i-1,i}b_1+b_{i-1,i}c_1+c_{i-1,i}a_1&=0\\
	a'_{ii}+b'_{ii}+c'_{ii}+a_1b_{i,i-1}+b_1c_{i,i-1}+c_1a_{i,i-1}&=0\\
	a_{i,i+1}+b_{i,i+1}+c_{i,i+1}+a_1b_{ii}+b_1c_{ii}+c_1a_{ii}&=0\\
	a_{i+1,i}+b_{i+1,i}+c_{i+1,i}+a'_{ii}b_1+b'_{ii}c_1+c'_{ii}a_1&=0
	\end{aligned}\right\}\text{ for } i\ge 1,
	\end{equation}
	where we set the convention $a_{01}=a_{10}=a_1, b_{01}=b_{10}=b_1, c_{01}=c_{10}=c_1$; and $a'_{11}=a_{11}, b'_{11}=b_{11},c'_{11}=c_{11}$, with the implication that the first and second equations coincide for $i=1$. Further, for the sake of ease in defining the automorphisms, we give meaningful interpretations for the equations in (\ref{eq:d3qDefnNew}) when  $i=0$. We set $a'_{00}=b'_{00}=c'_{00}=a_{00}=b_{00}=c_{00}=-1$; and $a_{0,-1}=b_{0,-1}=c_{0,-1}=a_{-1,0}=b_{-1,0}=c_{-1,0}=0$. Notice that the first and the second equations reduce to $-3=0$ for $i=0$, which is true in characteristic $3$.
	
	Now, we define five different linear maps on $\D_3(q)$ in Table \ref{table:autDq} below, by noting where each coordinate is sent to. For example, for fixed $x\in \bbF_q$, we denote $t_{1,1}(x)$ to be the map that satisfies $a_1\mapsto a_1+a_{-1,0}x = a_1$, $a_{11}\mapsto a_{11}+a_{00}x=a_{11}-x$, and so on. A ``-" as a table entry denotes a coordinate fixed by that map, e.g $t_{m+1,m}(a_{ii})=a_{ii}$.

	\begin{table}[ht]
		\renewcommand{\arraystretch}{1.2}
		\begin{tabularx}{\textwidth} { 
				| >{\centering\arraybackslash}X || >{\centering\arraybackslash}X | >{\centering\arraybackslash}X | >{\centering\arraybackslash}X |>{\centering\arraybackslash}X|>{\centering\arraybackslash}X|
			}
			\hline
			\small Coordinates $(i\ge 0)$& \small$t_{1,1}(x)$ & \small$t_{m,m+1}(x)$, $m\ge 1$; $r=i-m$ &\small$t_{m+1,m}(x)$, $m\ge 1$; $r=i-m$& \small$t_{m,m}(x)$, $m\ge 2$; $r=i-m$&\small$t'_{m,m}(x)$, $m\ge 2$; $r=i-m$\\ \hline \hline
			$a_{ii}$		&$+a_{i-1,i-1}x$ &$+a_{r,r-1}x$, $r\ge 1$ &- &$+a_{rr}x$, $r\ge 0$ &-\\	     \hline
            $a'_{ii}$		&$+a'_{i-1,i-1}x$ &- &$+a_{r-1,r}x$, $r\ge 1$ &- &$+a'_{rr}x$, $r\ge 0$\\    \hline
			$a_{i,i+1}$		&$+a_{i-1,i}x$ &$+a'_{rr}x$, $r\ge 0$ &- &$+a_{r,r+1}x$, $r\ge 0$ &-\\	     \hline
			$a_{i+1,i}$		&$+a_{i,i-1}x$ &- &$+a_{rr}x$, $r\ge 0$ &- &$+a_{r+1,r}x$, $r\ge 0$\\	     \hline
			 \hline
			$b_{ii}$		&$+b_{i-1,i-1}x$ &$+b_{r,r-1}x$, $r\ge 1$ &- &$+b_{rr}x$, $r\ge 0$ &-\\	     \hline
            $b'_{ii}$		&$+b'_{i-1,i-1}x$ &- &$+b_{r-1,r}x$, $r\ge 1$ &- &$+b'_{rr}x$, $r\ge 0$\\	 \hline
			$b_{i,i+1}$		&$+b_{i,i-1}x$ &$+b'_{rr}x$, $r\ge 0$ &- &$+b_{r,r+1}x$, $r\ge 0$ &-\\	     \hline
			$b_{i+1,i}$		&$+b_{i,i-1}x$ &- &$+b_{rr}x$, $r\ge 0$ &- &$+b_{r+1,r}x$, $r\ge 0$\\	     \hline
			 \hline
			$c_{ii}$		&$+c_{i-1,i-1}x$ &$+c_{r,r-1}x$, $r\ge 1$ &- &$+c_{rr}x$, $r\ge 0$ &-\\	     \hline
            $c'_{ii}$		&$+c'_{i-1,i-1}x$ &- &$+c_{r-1,r}x$, $r\ge 1$ &- &$+c'_{rr}x$, $r\ge 0$\\	 \hline
			$c_{i,i+1}$		&$+c_{i-1,i}x$ &$+c'_{rr}x$, $r\ge 0$ &- &$+c_{r,r+1}x$, $r\ge 0$ &-\\	     \hline
			$c_{i+1,i}$		&$+c_{i,i-1}x$ &- &$+c_{rr}x$, $r\ge 0$ &- &$+c_{r+1,r}x$, $r\ge 0$\\	     \hline
			
		\end{tabularx}
		\centering\caption{Automorphisms of $\D_3(q)$\label{table:autDq}}
		{\small($a'_{00}=b'_{00}=c'_{00}=a_{00}=b_{00}=c_{00}=-1$, $a_{0,-1}=b_{0,-1}=c_{0,-1}=a_{-1,0}=b_{-1,0}=c_{-1,0}=0$)}
	\end{table}	

    According to this convention, when $i=0$, the first two rows describe the images of $a_{00}$ and $a'_{00}$, which are constants and hence fixed by every map. The third and fourth rows coincide, and describe the images of $a_{1}$. When $i=1$, the first two rows coincide and describe the image of $a_{11}$. All other rows of Table~\ref{table:autDq} describe the images of unique coordinates.
    
	\begin{clm}
		\label{clm:autDq}
		The maps defined in Table \ref{table:autDq} are automorphisms of $\D_3(q)$.
	\end{clm}
	
	{\noindent \it Proof of Claim \ref{clm:autDq}.}
	First, we observe that each of the maps defined in Table \ref{table:autDq} is invertible when restricted to one of the vertex subsets $A$, $B$ or $C$.
    This is because when written as (infinite) matrices in the basis given by the coordinates, each map has $1$ along the diagonals and are lower triangular, thus are invertible.
    As an example, consider the action of $t_{1,1}$ on the vertex set $A$.
    When we write the matrix of $t_{1,1}$ in the standard basis, we obtain the following infinite lower-triangular matrix (here unfilled entries are $0$'s):
    \[
    t_{1,1} = \begin{bmatrix}
		1 &   &   &   &   &   &   & \\
		  0 & 1 &   &   &   &   &   & \\
		  0 & 0 & 1 &   &   &   &   & \\
		  0 & 0 & 0 & 1 &   &   &   & \cdots\\
		x & 0 & 0 &   & 1 &   &   & \\
		   & x & 0 &   &   & 1 &   & \\
		   &   & x &   &   &   & 1 & \\
        \vdots &  &  & \ddots &  &  & & \ddots\\
	\end{bmatrix}
    \]
    Thus, $t_{1,1}(x)$ is invertible for every $x\in\bbF_q$.
    A similar argument shows that all the maps in Table~\ref{table:autDq} are invertible.
    
	Hence, it remains to check that they are homomorphisms.
    Now, to show that a map $f$ is a homomorphism, it suffices to check that $\mathcal R^{(3)}(a,b,c) \implies \mathcal R^{(3)}(f(a), f(b), f(c))$, i.e. each relation in (\ref{eq:d3qDefnNew}) is preserved under $f$.
    We verify this implication for each map of Table~\ref{table:autDq} as follows.
	\begin{itemize}
		\item $t_{1,1}(x)$: We observe that the map $t_{1,1}(x)$ keeps $a_1,b_1,c_1$ fixed as $a_{1}=a_{0,1}\mapsto a_{0,1}+a_{-1,0}x = a_{0,1}$, etc. And, for $i\ge 1$, we need to check that the equations (\ref{eq:d3qDefnNew}) are preserved after the transformation given by $t_{1,1}$. Suppose the equations (\ref{eq:d3qDefnNew}) hold, then note that we also have for $i\ge 1$,
		\[
		\begin{array}{c}
		a_{ii}+b_{ii}+c_{ii}+a_{i-1,i}b_1+b_{i-1,i}c_1+c_{i-1,i}a_1=0,\\
		(a_{i-1,i-1}+b_{i-1,i-1}+c_{i-1,i-1}+a_{i-2,i-1}b_1+b_{i-2,i-1}c_1+c_{i-2,i-1}a_1)x=0,
		\end{array}
		\]
		and adding these up verifies that the first equation is preserved under the image of $t_{1,1}(x)$. Similarly, the other three equations can be verified for each $i\ge 1$.
		
		\item $t_{m,m+1}(x), m\ge 1$: Again, note that this map fixes $a_1=a_{0,1}$, $b_1=b_{0,1}$ and $c_1=c_{0,1}$ as for $i=0$ and $m\ge 1$, $r=i-m<0$. It also fixes all $a_{ii}$, $i\le m$ and all $a_{i,i+1}$, $i<m$. Therefore, all of (\ref{eq:d3qDefnNew}) are satisfied for $i<m$. When $i=m$, the first equation is still preserved as $a_{mm}, a'_{m-1,m}$ are fixed. For the third equation, we observe that $a_{m,m+1}\mapsto a_{m,m+1}+a'_{00}x=a_{m,m+1}-x$, $b_{m,m+1}\mapsto b_{m,m+1}-x$ and $c_{m,m+1}\mapsto c_{m,m+1}-x$. Thus, the third equation becomes
		\[
		(a_{m,m+1}-x)+(b_{m,m+1}-x)+(c_{m,m+1}-x)+a_1b_{mm}+b_1c_{mm}+c_1a_{mm}=0,
		\]
		which is still true as $3x=0$ in $\bbF_q$. Finally, for $i>m$, we need to check the validity of the first and third equations from (\ref{eq:d3qDefnNew}). However, note that for $i>m$ and $r=i-m\ge 1$,
		\[
		\begin{array}{c}
		a_{ii}+b_{ii}+c_{ii}+a_{i-1,i}b_1+b_{i-1,i}c_1+c_{i-1,i}a_1=0,\\
		(a_{r,r-1}+b_{r,r-1}+c_{r,r-1}+a'_{r-1,r-1}b_1+b'_{r-1,r-1}c_1+c'_{r-1,r-1}a_1)x=0,
		\end{array}
		\]
		and adding these up verifies the first equation, since $t_{m,m+1}(x)(a_{i-1,i})=a_{i-1,i}+a'_{r-1,r-1}x$. In a similar fashion, we verify the third equation by adding up:
		\[
		\begin{array}{c}
		a_{i,i+1}+b_{i,i+1}+c_{i,i+1}+a_1b_{ii}+b_1c_{ii}+c_1a_{ii}=0,\\
		(a'_{rr}+b'_{rr}+c'_{rr}+a_1b_{r,r-1}+b_1c_{r,r-1}+c_1a_{r,r-1})x=0,
		\end{array}
		\]
		for $i>m$ and $r=i-m\ge 1$. The second and fourth equations are unchanged by $t_{m,m+1}$.
		
		\item $t_{m+1,m}(x), m\ge 1$: Similar to $t_{m,m+1}$, this map fixes $a_{ii}$ and $a_{i,i+1}$ for every $i$, and hence does not change the first and third set of equations of (\ref{eq:d3qDefnNew}).
        For $i<m$, we have $r=i-m<0$, hence the map fixes all coordinates with $i<m$.
        For $i=m$, note that it changes $a_{m+1,m}\mapsto a_{m+1,m}-x$, yet fixes $a'_{mm}$.
        So, the second equation remains unchanged, and we also have
		\[
		(a_{m+1,m}-x)+(b_{m+1,m}-x)+(c_{m+1,m}-x)+a'_{mm}b_1+b'_{mm}c_1+c'_{mm}a_1=0.
		\]
        This shows that the fourth equation is preserved by the map.
		
		Finally, when $i>m$, the following four equations vouch for the validity of the second and fourth equations of (\ref{eq:d3qDefnNew}):
		\[
		\begin{array}{c}		
		\left\|
		\begin{array}{c}
		a_{i+1,i}+b_{i+1,i}+c_{i+1,i}+a'_{ii}b_1+b'_{ii}c_1+c'_{ii}a_1 = 0\\
		(a_{rr}+b_{rr}+c_{rr}+a_{r-1,r}b_1+b_{r-1,r}c_1+c_{r-1,r}a_1)x = 0\\
		\end{array}\right\|\\
		\left\|
		\begin{array}{c}
		a'_{ii}+b'_{ii}+c'_{ii}+a_1b_{i,i-1}+b_1c_{i,i-1}+c_1a_{i,i-1} = 0\\
		(a_{r-1,r}+b_{r-1,r}+c_{r-1,r}+a_1b_{r-1,r-1}+b_1c_{r-1,r}+c_1a_{r-1,r})x = 0\\
		\end{array}\right\|.
		\end{array}
		\]
		
		\item $t_{m,m}(x), m\ge 2:$ Same as before, we start by observing that $t_{m,m}(a_{mm})=a_{mm}-x$, $t_{m,m}(a_{m-1,m})=a_{m-1,m}$, preserving the first equation of (\ref{eq:d3qDefnNew}) for $i=m$. On the other hand, as $a_{m,m+1}\mapsto a_{m,m+1}+a_{0,1}x=a_{m,m+1}+a_1x$, we can rewrite the third equation into:
		\[
		(a_{m,m+1}+a_1x)+(b_{m,m+1}+b_1x)+(c_{m,m+1}+c_1x) + a_1(b_{mm}-x)+b_1(c_{mm}-x)+c_1(a_{mm}-x) = 0.
		\]
		
		For $i>m$ and $r=i-m\ge 1$, we only add the first and third equations to themselves for $i=i$ and $i=r$, after multiplying the $i=r$ equations by $x$.
		
		\item $t'_{m,m}(x),m\ge 2:$ For this map, $t'_{m,m}(a'_{mm}) = a'_{mm}-x$, $t'_{m,m}(a_{m,m-1})=a_{m,m-1}$, verifying the second equation of (\ref{eq:d3qDefnNew}) for $i=m$. And, as $t'_{m,m}(a_{m+1,m})=a_{m+1,m}+a_{1,0}x = a_{m+1,m}+a_1x$, we again have
		\[
		(a_{m+1,m}+a_1x)+(b_{m+1,m}+b_1x)+(c_{m+1,m}+c_1x) + (a'_{mm}-x)b_1+(b'_{mm}-x)c_1+(c'_{mm}-x)a_1 = 0.
		\]
		
		For $i>m$ and $r=i-m\ge 1$, adding the first and third equations to themselves for $i=i$ and $i=r$ completes the verification.
	\end{itemize}
	This calculation shows that the maps defined in Table \ref{table:autDq} are all homomorphisms.
    Since all of them are invertible, this finishes the proof of Claim \ref{clm:autDq}.\hfill{$\blacksquare$}
    
    \begin{rmk}
    It is worth mentioning here that the inverse of all the maps $f(x)$ are actually $f(-x)$ for $f\in\{t_{m,m+1}, t_{m+1,m}, t_{m,m}, t'_{m,m}\}$. 
    These matrices which have equal entries along the main diagonals are called Toeplitz matrices, and in general, their inverses may or may not be Toeplitz.
    For example, $t_{1,1}(x)\circ t_{1,1}(-x)$ is not the identity map.
    While the inverse of $t_{1,1}(x)$ when restricted to finite orders can be determined using, for example, \cite{lvToeplitz2007} or by induction, for our proof it suffices to just observe that $t_{1,1}(x)$ is invertible for $x\in\bbF_q$.
    \end{rmk}

	We now return to the proof of Proposition \ref{prop:a-to-a1-000}. In the proof of Claim \ref{clm:autDq}, we checked that $t_{1,1}(x)$ keeps $a_1$ fixed, and moves $a_{11}\mapsto a_{11}+a_{00}x = a_{11}-x$. Therefore, given any vertex $a\in A$ of $\D_3(q)$, we can perform $t_{1,1}(a_{11})$ to map $a_{11}$ to $0$. Let $a^{(11)} = t_{1,1}(a_{11})(a)$. Now, an application of $t_{1,2}(a^{(11)}_{12})$ sends the third coordinate, $a^{(11)}_{12}$ to $0$. Let $a^{(12)}=t_{1,2}(a^{(11)}_{12})(a^{(11)})$, and $a^{(21)}$, $a^{(22')}$, $a^{(32)},\ldots$ be defined similarly. Then, the map $\varphi$ given by
	\[
	\varphi = \cdots\circ t_{i+1,i}(a^{(i+1,i)}_{i+1,i})\circ t_{i,i+1}(a^{(ii')}_{i+1,i})\circ t'_{ii}(a^{(ii)\prime}_{ii})\circ t_{ii}(a^{(i-1,i)}_{ii})\circ \cdots \circ t_{1,2}(a^{(11)}_{12})\circ t_{1,1}(a_{11}),
	\]
	where $\varphi$ is truncated to $s$ compositions, sends the second through $(s+1)$-st coordinates of $a$ to $0$. It also preserves all edges through $a$, being an automorphism of $\D_3(q)$. This completes the proof. \hfill{$\Box$}
	

\subsubsection{Proof of Proposition~\ref{prop:D3kq-twolinks}}
Our goal in this section is to prove that $\D_3(k,q)$ admits at most two different link graphs.
    By Proposition~\ref{prop:a-to-a1-000}, it suffices to consider the link graphs of $a = (a_1,0,\ldots, 0)_A$ for $a_1\in\bbF_q$.
    Let $L_a$ denote the link graph of $a$.
    We see that $bc\in E(L_a)$ if and only if $\mathcal R_{k-1}(a,b,c)$ holds. This implies that the following equations hold ($i\ge 2$):
	\begin{equation}
	\label{eq:d3qLinkDefn}
    \mathcal R_{k-1}(a,b,c) = \left\{
	\begin{array}{ccc}
	b_{11}+c_{11}+a_1b_1+b_1c_1+c_1a_1&=&0\\
	b_{12}+c_{12}+a_1b_{11}+b_1c_{11}&=&0\\
	b_{21}+c_{21}+b_{11}c_1+c_{11}a_1&=&0\\
	&\vdots&\\
	b_{ii}+c_{ii}+b_{i-1,i}c_1+c_{i-1,i}a_1&=&0\\
	b'_{ii}+c'_{ii}+a_1b_{i,i-1}+b_1c_{i,i-1}&=&0\\
	b_{i,i+1}+c_{i,i+1}+a_1b_{ii}+b_1c_{ii}&=&0\\
	b_{i+1,i}+c_{i+1,i}+b'_{ii}c_1+c'_{ii}a_1&=&0\\
    &\vdots&
	\end{array}
    \right. 
	\end{equation}
	Here we consider two different cases.
	\begin{itemize}
		\item {\bf Case 1:} $a_1=0$.
		In this case, we note that the relations $\mathcal R_{k-1}(\vec 0_A, b, c)$ of (\ref{eq:d3qLinkDefn}) reduce to the relations $\mathcal R_{k-1}(b,c)$ of (\ref{eq:defDq-ab}) defining $\D(k,q)$, implying $L_a\cong \D(k,q)$.
	    
		\item {\bf Case 2:} $a_1\neq 0$.
		In this case, let us define an isomorphism $\psi:L_a\to L_{(1,0,\ldots,0)}$ as follows: 
		\[
        \arraycolsep=1.4pt\def\arraystretch{2.0}
	    \left\{\begin{array}{rl}
	    \psi(b_1)&=\dfrac{b_1}{a_1}\\
	    \psi(b_{ii})&=\dfrac{b_{ii}}{a_1^{2i}},\\ \psi(b'_{ii})&=\dfrac{b'_{ii}}{a_1^{2i}},\\ \psi(b_{i,i+1})&=\dfrac{b_{i,i+1}}{a_1^{2i+1}},\\ \psi(b_{i+1,i})&=\dfrac{b_{i+1,i}}{a_1^{2i+1}};
	    \end{array}\right\}
	    \mbox{ and }
	    \left\{\begin{array}{rl}
	    \psi(c_1)&=\dfrac{c_1}{a_1},\\
	    \psi(c_{ii})&=\dfrac{c_{ii}}{a_1^{2i}},\\ \psi(c'_{ii})&=\dfrac{c'_{ii}}{a_1^{2i}},\\ \psi(c_{i,i+1})&=\dfrac{c_{i,i+1}}{a_1^{2i+1}},\\ \psi(c_{i+1,i})&=\dfrac{c_{i+1,i}}{a_1^{2i+1}}.
	    \end{array}\right\}
		\]
		By dividing the equations in (\ref{eq:d3qLinkDefn}) by appropriate powers of $a_1$, it can be seen that $\psi$ is a homomorphism. As $a_1\neq 0$, $\psi$ is invertible, and hence $L_a\cong L_{(1,0,\ldots,0)}$, completing the proof. \hfill{$\Box$}
	\end{itemize}

\medskip
Proposition~\ref{prop:D3kq-twolinks} naturally leads us to investigate the links of the vertex $(1,0,\ldots 0)$ in $\D_3(q)$. Recall that the links of $(1,0,\ldots 0)_A$, $(1,0,\ldots 0)_B$ and $(1,0,\ldots 0)_C$ are isomorphic, so we may now consider $L_c$ where $c=(1,0,\ldots, 0)_C$.
The defining equations for this link is given by,
\begin{equation}
    \label{eq:d3qLinkDefnReduced}
    \mathcal R_{k-1}(a,b,c) = \left\{
    \begin{array}{ccc}
    a_{11}+b_{11}+a_1+a_1b_1+b_1&=&0\\
    a_{12}+b_{12}+a_{11}+a_1b_{11}&=&0\\
    a_{21}+b_{21}+a_{11}b_1+b_{11}&=&0\\
    &\vdots&\\
    a_{ii}+b_{ii}+a_{i-1,i}b_1+b_{i-1,i}&=&0\\
    a'_{ii}+b'_{ii}+a_{i,i-1}+a_1b_{i,i-1}&=&0\\
    a_{i,i+1}+b_{i,i+1}+a_{ii}+a_1b_{ii}&=&0\\
    a_{i+1,i}+b_{i+1,i}+a'_{ii}b_1+b'_{ii}&=&0\\
    &\vdots&
    \end{array}\right.
\end{equation}

We can reduce this further by replacing $a_1$ with $a_1+1$ and $b_1$ with $b_1+1$. Noting that $(a_1+1)+(a_1+1)(b_1+1)+(b_1+1) = a_1b_1 -a_1 -b_1$ in characteristic $3$, we get a new set of equations, namely (\ref{eq:defDpq-ab}).
We call this new series of graphs $\D'(k,q)$, and take a closer look at them in the next subsection.

\subsection{The bipartite graphs $\D'(k,q)$}
We now take a detour into the sequence of graphs $\D'(k,q)$.
It is worth clarifying that in this subsection, we look at $\bbF_q$ of arbitrary finite characteristic.
\begin{defn}[The bipartite graphs $\D'(q)$]
	For a prime power $q$, let $A$ and $B$ be two \textit{disjoint} copies of the countably infinite dimensional vector space $V$ over $\bbF_q$. We use the following coordinate representations for $a\in A$, $b\in B$:
	\begin{equation}
	\label{eq:coordinatesDpq}
	\begin{array}{lllcclcll}
	a&=&(&\overbrace{a_1,a_{11},a_{12},a_{21}}^{\text{group }1}, &\overbrace{a_{22},a'_{22},a_{23},a_{32}}^{\text{group }2},&\ldots&\overbrace{a_{ii},a'_{ii},a_{i,i+1},a_{i+1,i}}^{\text{group }i},&\ldots&),\\
	b&=&(&b_1,b_{11},b_{12}, b_{21}, &{b_{22},b'_{22},b_{23},b_{32}},&\ldots &b_{ii},b'_{ii},b_{i,i+1},b_{i+1,i},&\ldots&).
	\end{array}
	\end{equation}
	Let $\D'(q)$ consist of vertex set $A\sqcup B$, and let us join $a\in A$ to $b\in B$ iff the following equations hold ($i\ge 2$):
	\begin{equation}
	\label{eq:defDpq-ab}
    \mathcal R'(a,b):=\left\{
	\begin{array}{ccc}
	a_{11}-a_1+b_{11}-b_1+a_1b_1&=&0\\
	a_{12}+a_{11}+b_{12}+b_{11}+a_1b_{11}&=&0\\
	a_{21}+a_{11}+b_{21}+b_{11}+a_{11}b_1&=&0\\
	&\vdots&\\
	a_{ii}+a_{i-1,i}+b_{ii}+b_{i-1,i}+a_{i-1,i}b_1&=&0\\
	a'_{ii}+a_{i,i-1}+b'_{ii}+b_{i,i-1}+a_1b_{i,i-1}&=&0\\
	a_{i,i+1}+a_{ii}+b_{i,i+1}+b_{ii}+a_1b_{ii}&=&0\\
	a_{i+1,i}+a'_{ii}+b_{i+1,i}+b'_{ii}+a'_{ii}b_1&=&0\\
    &\vdots&
	\end{array}\right.
	\end{equation}

    Again, we observe that the first $k-1$ equations of (\ref{eq:defDpq-ab}) has the first $k$ coordinates of (\ref{eq:coordinatesDpq}).
    So if $A_k$ and $B_k$ denote the truncation of $A$ and $B$ from (\ref{eq:coordinatesDq}) to the first $k$ coordinates, and $\mathcal R'_k$ the truncation of $\mathcal R'$ from (\ref{eq:defDq-ab}) to the first $k$ relations, then, we define $\D'(k,q)$ as the bipartite graph with bipartition $A_k\sqcup B_k$ where vertices $a\in A_k$ and $b\in B_k$ are adjacent iff $\mathcal R'_{k-1}(a,b)$.

    By an exactly analogous argument as for $\D(k,q)$, it follows that $\D'(k,q)$ is a $q$-regular graph on $2q^k$ vertices.
\end{defn}

It is natural to inquire whether $\D'(k,q)$ and $\D(k,q)$ are related in any way, in particular, whether they're the same graph.
The answer turns out to be yes for small values of $k$, but no for larger $k$:
\begin{thm}
	\label{thm:Dpkq-Dkq-isomorphism}\hfill
 
	(a) For $2\le k\le 6$, $\D'(k,q)\cong \D(k,q)$.
	
	(b) $\D'(11,3)\not\cong \D(11,3)$.
\end{thm}
\begin{proof}
	First, we prove part (a).
	
	The main idea of the proof is as follows.
	Observe that it is enough to show that $\D'(6,q)\cong \D(6,q)$, as an isomorphism $\D'(6,q)\to \D(6,q)$ can be restricted to fewer coordinates to give isomorphisms $\D'(k,q)\to \D(k,q)$ for $k\le 6$.
	To demonstrate that $\D'(6,q)\cong \D(6,q)$, we shall define a map $x\mapsto \overbar x$ sending $x\in V(\D'(6,q))$ to the vector $\overbar x\in \bbF_q^{6}$, such that for $a\in A$ and $b\in B$, we have $ab\in E(\D'(6,q))$ implies $\overbar a \overbar b \in E(\D(6,q))$.
	By construction, this map will be linear and invertible, which would then complete the proof.
	
	We define the map $x\mapsto \overbar x$ as described in Table~\ref{table:Dpkq-Dkq-isomorphism}.
	
	\begin{table}[ht]
		\renewcommand{\arraystretch}{1.2}
		\begin{tabularx}{\textwidth} { 
				| >{\centering\arraybackslash}X | >{\centering\arraybackslash}X || >{\centering\arraybackslash}X | >{\centering\arraybackslash}X|
			}
			\hline
			$a\in V(\D'(6,q))\cap A$ & $\overbar a\in \bbF_q^{10}$ & $b\in V(\D'(6,q))\cap B$ & $\overbar b\in \bbF_q^{10}$\\ \hline \hline
			$a_1$ & $a_1$ & $b_1$ & $b_1$\\ \hline
			$a_{11}$ & $a_{11}-a_1$ & $b_{11}$ & $b_{11}-b_1$\\ \hline
			$a_{12}$ & $a_{12}+a_1$ & $b_{12}$ & $b_{12}+b_1$\\ \hline
			$a_{21}$ & $a_{21}+a_1$ & $b_{21}$ & $b_{21}+b_1$\\ \hline
			$a_{22}$ & $a_{22}+a_{12}+a_{11}-a_1$ & $b_{22}$ & $b_{22}+b_{12}+b_{11}-b_1$\\ \hline
			$a'_{22}$ & $a'_{22}+a_{21}+a_{11}-a_1$ & $b'_{22}$ & $b'_{22}+b_{21}+b_{11}-b_1$\\ \hline
		\end{tabularx}
		\centering\caption{The isomorphism $\D'(6,q)\to \D(6,q)$}
		\label{table:Dpkq-Dkq-isomorphism}
	\end{table}

	Suppose $a,b\in V(\D'(k,q))$ with $a\in A$, $b\in B$ and $ab\in E(\D'(k,q))$. This implies:
	\[
	\begin{array}{ccc}
	a_{11}-a_1+b_{11}-b_1+a_1b_1 &=& 0\\
	a_{12}+a_{11}+b_{12}+b_{11}+a_1b_{11}&=&0\\
	a_{21}+a_{11}+b_{21}+b_{11}+a_{11}b_1&=&0\\
	a_{22}+a_{12}+b_{22}+b_{12}+a_{12}b_1&=&0\\
	a'_{22}+a_{21}+b'_{22}+b_{21}+a_1b_{21}&=&0
	\end{array}
	\]
	Now observe that, $\overbar a_1 = a_1$ and $\overbar b_1 = b_1$. Further,
    
    \vspace{-20pt}
	\begin{equation}
    \arraycolsep=1.0pt\def\arraystretch{4.0}
	\begin{array}{rl}
		\bullet & \left\{\begin{aligned}\overbar a_{11}+\overbar b_{11}+a_1b_1 &= a_{11}-a_1+b_{11}-b_1 +a_1b_1 \\& = 0,\end{aligned} \right.\\
		\bullet & \left\{\begin{aligned}
		\overbar a_{12}+\overbar b_{12}+a_1\overbar b_{11} &= a_{12}+a_1+b_{12}+b_1+a_1(b_{11}-b_1)\\
		&=a_{12}+a_1 + b_{12}+b_1 + a_1b_{11} + (a_{11}-a_1 + b_{11}-b_1)\\
		&=a_{12}+a_{11}+b_{12}+b_{11}+a_1b_{11} \\ &= 0,
		\end{aligned}\right. \\
		\bullet & \left\{\begin{aligned}
		\overbar a_{21} + \overbar b_{21} + \overbar a_{11}b_1 &= a_{21}+a_1 + b_{21}+b_1 + (a_{11}-a_1)b_1 \\
		& = a_{21}+a_1 + b_{21}+b_1 +a_{11}b_1 + (a_{11}-a_1+b_{11}-b_1)\\
		& = a_{21}+a_{11} + b_{21} + b_{11} + a_{11}b_1\\
		& = 0, 
		\end{aligned}\right. \\
		\bullet & \left\{\begin{aligned}
		\overbar a_{22} + \overbar b_{22} + \overbar a_{12}b_1 & = a_{22}+a_{12}+a_{11}-a_1  + b_{22} + b_{12} + b_{11}-b_1 + (a_{12}+a_1)b_1\\
		&= a_{22} + a_{12} + b_{22} + b_{12} +a_{12}b_1\\
		&=0,
		\end{aligned}
		\right. \\
		\bullet & \left\{\begin{aligned}
		\overbar a'_{22} + \overbar b'_{22} + \overbar a_{1}\overbar b_{21} & = a'_{22}+a_{21}+a_{11}-a_1  + b'_{22} + b_{21} + b_{11}-b_1 + a_1(a_{21}+b_1)\\
		&= a'_{22} + a_{21} + b'_{22} + b_{21} +a_{1}b_{21}\\
		&=0.
		\end{aligned}
		\right.
	\end{array}
	\end{equation}
	Therefore the map $x\mapsto \overbar x$ is an isomorphism from $\D'(6,q)$ to $\D(6,q)$, as desired. \hfill{$\blacksquare$}
	
	Our proof of part (b) is purely computational. In summary, it has been computed that the diameter of the component of $\D(11,3)$ containing $\vec 0$ is $22$ whereas the same number for $\D'(11,3)$ is $20$, implying they're not isomorphic (as it is known that $\D(11,3)$ is edge-transitive). Further, $\D(11,3)$ has $112$ cycles through the edge $\{\vec 0, \vec 0\}$ whereas $\D'(11,3)$ has only $4$. This also implies $\D(11,3)\not\cong \D'(11,3)$.
	
	The github repository \url{https://github.com/Potla1995/hypergraphSuspension/} contains further details on how to reproduce these results.
\end{proof}

\begin{rmk}
	Computer calculations for small values of $q$ suggest that $\D'(k,q)$ and $\D(k,q)$ are isomorphic for $7\le k\le 10$. However, the proof method used for $k\le 6$ does not extend to this range.
\end{rmk}

Note that proving that $\D'(k,q)$ has high girth is synonymous to proving lower bounds on $\ex(n,\widetilde C_{2k})$ by the machinery we've built so far in this section.
There is computational evidence for up to $k=13$ that the girth of $\D'(k,q)$ is at least $k+4$ if $k$ is even, and $k+5$ if $k$ is odd, analogous to $\D(k,q)$.
As $\D'(k,q)$ is a sequence of graphs not isomorphic with $\D(k,q)$ in general, we propose to study the following open question:

\begin{prob}
	\label{prob:girthdpkq}
	What is the girth of $\D'(k,q)$ and connectivity of $\D'(k,q)$ for values of $k\ge 7$?
\end{prob}

\subsection{Proof of Theorem~\ref{thm:c2ktilde}}
We have now built all the machinery required to complete our proof of Theorem~\ref{thm:c2ktilde}, and will delve into the proof.

\begin{proof}[Proof of Theorem \ref{thm:c2ktilde}]
Recall that we have to check three properties of $\D_3(k,q)$, and that $q$ is a power of $3$.
	\begin{enumerate}
	    \item First, we check that $\D_3(k,q)$ has $3q^k$ vertices and $q^{2k+1}$ edges.
	    It is clear that every part of $\D_3(k,q)$ has $q^k$ vertices.
	    Since there is exactly one free variable when we fix $a$ and $b$ for a hyperedge $\{a,b,c\}$, this gives us a total of $q^k\cdot q^k\cdot q=q^{2k+1}$ edges.
	    
	    \item Next, we shall prove that the link graphs of every vertex of $\D_3(k,q)$ is isomorphic, in fact, to $\D(k,q)$ for $k\le 6$.
	    By Proposition~\ref{prop:D3kq-twolinks}, the link of every vertex of $\D_3(k,q)$ is isomorphic to $\D(k,q)$ or $\D'(k,q)$ as $q$ is a power of $3$.
	    However, $\D(k,q)\cong \D'(k,q)$ for $k\le 6$, implying the required assertion.
	    
	    \item Finally, it remains to show that $\D_3(3,q)$ is $\widetilde C_6$-free and $\D_3(5,q)$ is $\widetilde C_8$-free. From the previous point, and since $\D(3,q)$ and $\D(5,q)$ are known to have girths 8 and 10 respectively (Proposition~\ref{prop:dkqproperties} part 2), this completes the proof.
	\end{enumerate}
 \vspace{-27pt}
\end{proof}

\section{The arc construction and Wenger's construction}
In this section, we relate the arc construction and Wenger's construction via Proposition~\ref{prop:Arc_Wenger_isomorphism}, and provide a new set of $C_6$-free graphs with $n$ vertices and $\Theta(n^{4/3})$ edges via proving Theorem~\ref{thm:newC6construction}.

\subsection{Proof of Proposition~\ref{prop:Arc_Wenger_isomorphism}}
Our main goal is to algebraically parametrize the constructions $G_{\tt arc}(k,q,\alpha_0)$ for $k\ge 2$, prime powers $q$ and the normal rational curve $\alpha_0$, which would lead us to Wenger's construction $H(k,q)$.
To this end, we would require the use of the Pl\"ucker embedding~\cite{griffithsHarrisAG1978}, an algebraic geometric tool that allows us to parametrize the set $L$.

\begin{lemma}[Pl\"ucker Embedding]
\label{lem:plucker}
Every line $\ell$ passing through points $[a_1:\cdots:a_{t+1}]$ and $[b_1:\cdots:b_{t+1}]$ in $\PG(t,q)$ can be parametrized using $\binom {t+1}2$ coordinates $\{w_{ij}: 1\le i<j \le t+1\}$, where $w_{ij}$ is given by the determinant of the $2\times 2$ matrix obtained by appending the $i$'th and $j$'th rows of
\[
\begin{bmatrix}
a_1 & a_2 &\cdots & a_{t+1}\\
b_1 & b_2 &\cdots & b_{t+1}
\end{bmatrix},
\]
i.e. $w_{ij} = a_ib_j - a_jb_i$.
\end{lemma}

For further details on the Pl\"ucker embedding, the reader is referred to \cite{griffithsHarrisAG1978}, p.211.

We are now well-equipped to prove Proposition~\ref{prop:Arc_Wenger_isomorphism}, which asserts that $G_{\tt arc}(k,q,\alpha_0^-)\cong H(k,q)$.

\begin{proof}[Proof of Proposition~\ref{prop:Arc_Wenger_isomorphism}]
Recall that in the $G_{\tt arc}(k,q,\alpha_0^-)$ construction, $P=\Sigma\setminus\Sigma_0$ and \[L=\{\mbox{projective lines }\ell: \ell\cap\Sigma_0 \in \alpha_0^-\}.\]
Therefore, $|P|=q^k$ and $|L|=q^{k-1}|\alpha_0^-| = q^k$.

Observe that any line in $L$ passes through a point $[1:a_1:\cdots:a_{k}]\in P$ and a point $[0:1:x:\cdots:x^{k-1}]\in\alpha_0^-$. Let $\{w_{ij}:1\le i<j\le k+1\}$ parametrize lines in $L$. Then, for $2\le j\le k+1$,
\begin{equation}
\label{eq:pluckerGrassmannian_w1j}
w_{1j} = \det \begin{bmatrix}1&a_{j-1}\\0&x^{j-2}\end{bmatrix} = x^{j-2},
\end{equation}
and for $2\le i < j$,
\begin{equation}
\label{eq:pluckerGrassmannian_wij}
w_{ij}=\det\begin{bmatrix}a_{i-1}&a_{j-1}\\x^{i-2}&x^{j-2}\end{bmatrix}=a_{i-1}x^{j-2}-a_{j-1}x^{i-2}.
\end{equation}
Equation (\ref{eq:pluckerGrassmannian_w1j}) implies that $w_{13}=x$ and $w_{1j}=w_{13}^{j-2}$ for $j\ge 2$.
Moreover, plugging $i=2$ into (\ref{eq:pluckerGrassmannian_wij}) gives us $w_{2j} = a_1x^{j-2}-a_{j-1}$ for $j>2$.
Thus, $a_{j-1} = a_1w_{13}^{j-2}-w_{2j}$ for $j\ge 3$.
Now, for any $3\le i < j$, we have
\[
w_{ij} = (a_1w_{13}^{i-2}-w_{2i})w_{13}^{j-2} - (a_1w_{13}^{j-2}-w_{2j})w_{13}^{i-2}=w_{13}^{i-2}w_{2j} - w_{13}^{j-2}w_{2i}.
\]

In particular, the above analysis implies that $w_{1j}$ are all dependent on $w_{13}$ and $\{w_{ij}: i\ge 3\}$ are all dependent on $w_{13}$ and $\{w_{2j}:j\ge 3\}$.
Hence we may reduce our free variables to only the set $\{w_{13}\}\cup \{w_{2j}:3\le j\le k+1\}$.
Let $b_1 := x = w_{13}$ and $b_{j-1} = w_{2j}, 3 \le j \le k+1$. Then, the equation (\ref{eq:pluckerGrassmannian_wij}) for $i=2$ reduces to
\[
b_{j-1} = a_1b_1^{j-2}-a_{j-1}, \ 3\le j \le k+1,
\]
Which is exactly the defining set of equations for the graph $H(k,q)$. As $P$ consists of $q^k$ points parametrized by $\{w_{13}\}\cup \{w_{2j}:3\le j\le k+1\}$, this implies $G_{\tt arc}(k,q,\alpha_0^-)\cong H(k,q)$.
\end{proof}

\subsection{Proof of Theorem~\ref{thm:newC6construction}}
We remark that Theorem~\ref{thm:newC6construction} can be proved completely analogously to the proof of Proposition~\ref{prop:Arc_Wenger_isomorphism} via using the arc $\alpha$ of $\PG(2,2^r)$ given by $\alpha = \{[1:t:t^{2^s}]: t\in \mathbb F_q\}$. 
However, for the sake of simplicity, 
we provide an alternative and more direct proof following Wenger's proof in~\cite{wengerExtremal1991}. Recall that $q=2^r$, $(s,r)=1$, and $G(2^r, s)$ is the bipartite graph with parts $A=B=\bbF_q^3$ such that $(a_1,a_2,a_3)\in A$ and $(b_1,b_2,b_3)\in B$ are adjacent iff
\[
b_2+a_2 = a_1b_1 \mbox{ and } b_3+a_3 = a_1b_1^{2^s}.
\]

\begin{proof}[Proof of Theorem~\ref{thm:newC6construction}]
	Let $a = (a_1,a_2,a_3), b = (b_1,b_2,b_3), \ldots, f = (f_1,f_2,f_3)$ form a $C_6$ in $G(2^r,s)$ where $a,c,e\in A$ are distinct, and $b,d,f\in B$ are distinct.
	
	Then, as $ab$ and $bc$ are edges, we have $a_2+b_2 = a_1b_1, c_2+b_2 = c_1b_1$ implying $a_2+c_2 = b_1(a_1+c_1)$ (due to characteristic $2$). Similarly, $a_3+c_3 = b_1^{2^s}(a_1+c_1)$. We can write these equations as,
	\[
	\begin{bmatrix}
	a_1+c_1\\
	a_2+c_2\\
	a_3+c_3
	\end{bmatrix}=
	\begin{bmatrix}
	1 \\
	b_1 \\
	b_1^{2^s}
	\end{bmatrix}\cdot (a_1+c_1),
	\]
	and similarly 
	\[
	\begin{bmatrix}
	c_1+e_1\\
	c_2+e_2\\
	c_3+e_3
	\end{bmatrix}=
	\begin{bmatrix}
	1 \\
	d_1 \\
	d_1^{2^s}
	\end{bmatrix}\cdot (c_1+e_1)\mbox{ and }
	\begin{bmatrix}
	e_1+a_1\\
	e_2+a_2\\
	e_3+a_3
	\end{bmatrix}=
	\begin{bmatrix}
	1 \\
	f_1 \\
	f_1^{2^s}
	\end{bmatrix}\cdot (e_1+a_1).
	\]
	Adding these up and using characteristic 2, we have
	\begin{equation}
	\label{eq:matrixC6free}
	\begin{aligned}
	\begin{bmatrix}0\\0\\0\end{bmatrix}&=
	\begin{bmatrix}
	1 \\
	b_1 \\
	b_1^{2^s}
	\end{bmatrix}\cdot (a_1+c_1)+\begin{bmatrix}
	1 \\
	d_1 \\
	d_1^{2^s}
	\end{bmatrix}\cdot (c_1+e_1)+\begin{bmatrix}
	1 \\
	f_1 \\
	f_1^{2^s}
	\end{bmatrix}\cdot (e_1+a_1)\\
	& = \begin{bmatrix} 1& 1 & 1\\ b_1 & d_1 & f_1 \\ b_1^{2^s} & d_1^{2^s} & f_1^{2^s}\end{bmatrix}\begin{bmatrix} a_1+c_1 \\ c_1+e_1\\ e_1+a_1\end{bmatrix}.
	\end{aligned}.
	\end{equation}
	Let $M(x,y,z) = { \begin{bmatrix} 1& 1 & 1\\ x & y & z \\ x^{2^s} & y^{2^s} & z^{2^s}\end{bmatrix}}$.
	
	\begin{clm}
		\label{clm:M(x,y,z)Invertible}
		If $x,y,z\in\bbF_q$ are all distinct, then $M(x,y,z)$ is invertible.
	\end{clm}
	{\it Proof of Claim~\ref{clm:M(x,y,z)Invertible}.}
	We wish to show that $\det M(x,y,z)\neq 0$, which simplifies to $\frac{x^{2^s}+y^{2^s}}{x+y}\neq \frac{y^{2^s}+z^{2^s}}{y+z}$. If there were pairwise distinct $t_1,t_2,t_3\in \bbF_q$ with $\frac{t_1^{2^s}+t_2^{2^s}}{t_1+t_2}=\frac{t_2^{2^s}+t_3^{2^s}}{t_2+t_3}$, we would then have, for a fixed $t_2$ and for $x=t_1+t_2$ and $y=t_3+t_2$, that $\frac{(x+t_2)^{2^s}+t_2^{2^s}}{x} = \frac{(y+t_2)^{2^s}+t_2^{2^s}}{y}$. This would imply that  $\left|\left\{\frac{(x+t_2)^{2^s}+t_2^{2^s}}{x}: x\in\bbF_q\setminus\{t_2\} \right\}\right| < q-1$. 
    Therefore, it is enough to check that for any arbitrary $t\in\bbF_q$,
	\[
	\left|\left\{\frac{(x+t)^{2^s}+t^{2^s}}{x}: x\in\bbF_q\setminus\{t\} \right\}\right| = q-1.
	\]
	Observe that, by the binomial theorem and using the fact that $\binom {2^s}i$ is even for every $0<i<2^s$, $\frac{(x+t)^{2^s}+t^{2^s}}{x} = x^{2^s-1}$. Hence, it suffices to show that the map $x\mapsto x^{2^s-1}$ is a permutation of $\bbF_q$. However, as the multiplicative group $\bbF_q^\ast$ has order $q-1$, this happens only when $(2^s-1,q-1)=1$, which is true since
	\[
	(2^s-1,2^r-1) = 2^{(s,r)}-1 = 1,
	\]
	by assumption.\footnote{Here we use the elementary fact that $(a^m-1,a^n-1)=a^{(m,n)}-1$ for any natural numbers $a$, $m$, $n$. This can be shown by iteratively using the euclidean algorithm: if $n\ge m$, $(a^m-1, a^n-1) = (a^m-1, a^n-a^m)=(a^m-1, a^{n-m}-1)$.} \hfill{$\blacksquare$}
	
	Now, we see that $b_1\neq d_1$. This is since if $b_1=d_1$, then, as 
	\[
	b_2+c_2=b_1c_1=c_1d_1=c_2+d_2
	\]
	and 
	\[
	b_3+c_3 = b_1^{2^s}c_1=c_1d_1^{2^s} = c_3+d_3,
	\]
	we would obtain $ b =  d$, a contradiction. Thus, $b_1, d_1, f_1$ are pairwise distinct, and therefore $M(b_1,d_1,f_1)$ is invertible. Hence, (\ref{eq:matrixC6free}) implies 
	\[
	a_1+c_1 = c_1+e_1 = e_1+a_1 = 0,
	\]
	i.e., $a_1=c_1=e_1$. However, as 
	\[
	a_2+b_2 = a_1b_1 = c_1b_1 = b_2+ c_2
	\]
	and 
	\[
	a_3+b_3 = a_1b_1^{2^s} = c_1b_1^{2^s} = b_3+c_3,
	\]
	this would imply $ a =  c$, a contradiction.
\end{proof}

\bigskip 

{\bf Acknowledgments.} This work was supported by the Center of Innovations for Sustainable Quantum AI (JST Grant Number JPMJPF2221). The author is immensely grateful to Dhruv Mubayi for several helpful discussions and warm encouragement, and Dhruv Mubayi, Xizhi Liu and an anonymous referee for several comments and suggestions that improved the presentation of the paper.

\printbibliography



%
%
%
%

\end{document}